\documentclass[doublsp,12]{algorithmica}

\begin{document}
\newcommand{\tab}[0]{\hspace{.1in}}
\newtheorem{thm}{\noindent {\sc Theorem}}
\newtheorem{lem}{\noindent {\sc Lemma}}
\newenvironment{proof}{\noindent {\sc Proof. }}
  {\hspace{8in minus8in}$\Box$\vspace{.1in}}

\pagestyle{myheadings}
\markboth{}{} 

\title{\vspace{-1in}Exact Sampling from Perfect Matchings of Dense 
Nearly Regular Bipartite Graphs\footnote{Supported by NSF postdoc 99-71064} }
\author{Mark Huber\footnote{Department of Mathematics and ISDS, 
Duke University, Box 90320, Durham, NC  27708-0320, {\tt mhuber@math.duke.edu} 
} } 
\date{}
\maketitle
\thispagestyle{empty}

%{\small \noindent {\bf Abstract}
\begin{abstract}
We present the first algorithm for generating random variates
exactly uniformly from the set of perfect matchings of a bipartite
graph with a polynomial expected running time over a nontrivial
set of graphs.  Previous Markov chain results obtain 
approximately uniform variates for arbitrary graphs in polynomial
time, but their general running time is $\Theta(n^{26} (\ln n)^2).$
Our algorithm employs acceptance/rejection together with a new upper
limit on the permanent of a form similar to Bregman's Theorem.  
For a graph with $2n$ nodes where
the degree of every node is nearly 
$\gamma n$ for a constant $\gamma$, the
expected running time is 
$O(n^{1.5 + .5/\gamma})$. 
Under these conditions, Jerrum and Sinclair showed that a Markov chain
of Broder can generate approximately uniform 
variates in $\Theta(n^{4.5 + .5/\gamma} \ln n)$ time,
making our algorithm significantly faster on this class of graph.
With our approach, approximately counting the number of perfect
matchings (equivalent to finding the permanent of a 0-1 matrix and so
$\sharp P$ complete) can
be done without use of selfreducibility.
We give an $1 + \sigma$ approximation 
algorithm for finding the permanent of 0-1
matrices with nearly 
identical row and column sums that runs in time
$O(n^{1.5 + .5/\gamma}
 \frac{1}{\sigma^2}
 \log \frac{1}{\delta})$, where the
probability the output is
within $1 + \sigma$ of the permanent is at least $1- \delta$.
%}
\end{abstract}


\section{The Problem.}
\label{SEC:I}
\pagenumbering{arabic}
Given a bipartite graph $(V,E)$ with $n$ 
nodes in each of the two partitions (so $|V| = 2n$),
a perfect matching is a subset of edges such that each node is adjacent
to exactly one edge.  Let $\Omega$ denote the set of all such perfect
matchings.  Generating random variates uniformly from $\Omega$ has a variety
of applications, such as determining the $p$ values of tests for
doubly truncated data~\cite{efronp1999}.  For general graphs,
Jerrum, Sinclair, and Vigoda~\cite{jerrumsv2001} gave the first polynomial
time algorithm for generating 
samples that are arbitrarily close to uniform, 
but their method requires
$\Theta(n^{26} (\ln n)^2)$ steps. 

Generation algorithms can also be used in constructing algorithms
for approximately counting $\Omega$.  This problem
is equivalent to finding the permanent of a 0-1 matrix.
Computing the permanent of a matrix was
one of the first problems shown to be $\sharp P$ 
complete~\cite{valiant79}.  
When restricted to 0-1 matrices 
where the row and column sums are at least 
$n / 2$, the problem remains $\sharp P$ complete.
We restrict further the problem to the case where each edge has at 
least $n \gamma$ edges and the graph is nearly regular.  Even under all
these restrictions, the problem remains $\sharp P$ complete~\cite{broder86}.

Algorithms for approximating the permanent have been
constructed using a variety of techniques, including decomposition of
the problem~\cite{jerrumv96} and using determinants of related random 
matrices~\cite{barvinok1999}.  Using selfreducibility of the 
problem~\cite{jerrumvv86}, any polynomial time technique for generating
samples can be used to efficiently approximate $|\Omega|$, but the number
of samples needed can be fairly large.  The Jerrum, Sinclair,
and Vigoda method constructs such an estimate along with a means
for generating samples simultaneously--the first sample requires
the most time, making the approximation algorithm $\Theta(n^{26}(\ln n)^2)$
as well.

For bipartite graphs
with minimum degree $n/2$, the Markov chain of Broder~\cite{broder86}
requires 
$O(n^5 \ln n)$ time to generate approximate 
samples~\cite{diaconiss91}, and will more generally be polynomial when
the ratio of the number of almost matchings (matchings with
$n-1$ edges) to the number of perfect matchings is 
polynomial~\cite{sinclair92}.

Given the $\Theta(n^{26} (\ln n)^2)$ technique for generating 
approximate samples, two questions arise.  First, can this
running time be significantly reduced by restricting the
class of graphs we consider?  Second, can samples be drawn exactly
from the uniform distribution, as opposed to the approximate samples
created by Markov chains?

An exact sampler is preferable to approximate samplers for several
reasons.  A desirable property of stochastic algorithms is that they be
unbiased, that the expected value of the output be equal to the
actual output.  Finding unbiased algorithms for applications such
as estimating $p$ values relies on the output coming exactly from
the desired distribution.  A more important reason in practice is that
algorithms coming from Markov chains are not only $O(T)$, they are
$\Theta(T)$ where $T$ is
the mixing time of the chain.  That is, $T$ steps must always be
taken in order to insure that the output is close to the desired
distribution.  On the other hand, exact sampling algorithms are
often truly $O(T')$, where $T'$ is merely an upper bound on the running time
that in practice is never reached.

In the remainder of this paper, we present an exact sampling
technique that is guaranteed to run in polynomial time over a 
nontrivial set of graphs.  Roughly speaking, when
the degree of each node is close to $\gamma n$, the algorithm
takes time 
$O(n^{1.5 + 0.5 / \gamma})$.

Also, this method naturally gives an algorithm for approximately
counting the number of perfect matchings without the need to take
advantage of selfreducibility.  This makes implementation easier, and
also makes the running time similar to that needed to draw a sample:
the running time is
$O(S \frac{1}{\epsilon^2}\log\left(\frac{1}{\delta}\right))$ time,
where $S$ is the expected time needed to obtain a single random variate. 

The permanent of a matrix $A$ is
defined as 
$$per(A) := \sum_{\pi \in S_n} \prod_{i = 1}^n a_{i \pi(i)}.$$
When the matrix is 0-1 the terms of the permanent are 0 or 1.
They are nonzero precisely when $\{\{1,\pi(1)\},\ldots,\{n,\pi(n)\}\}$ 
is a perfect matching in the bipartite graph where there is an edge from the 
$i$th node of one partition to the $j$th node of the other partition if and
only if $a_{ij} = 1$.
This makes computing the permanent of a 0-1 matrix
equivalent to
counting the number of perfect matchings.

We will only be able to draw from perfect matchings where
the graph is nearly regular. Although our method can
be run on instances which are far from being regular and will generate
an exact sample when it terminates, it is not guaranteed to terminate
in polynomial time unless the graph is nearly regular
with approximate degree $\Delta$ that is at least $\gamma n$.  

\subsection{Bregman's Theorem.}
The heart of our work is a new inequality for permanents that is of the
same form but slightly weaker than
Bregman's Theorem.  Although Bregman's theorem is stronger, it cannot
be used algorithmically to generate samples, whereas our 
new inequality can be used in such a fashion.

First conjectured by Minc~\cite{minc1963} and proved by 
Bregman~\cite{bregman73}, Bregman's Theorem
gives the following upper bound on the permanent:
\begin{equation}
\label{EQN:mincs}
M(A) := \prod_{i=1}^n (r(i)!)^{1/r(i)} \geq per(A),
\end{equation}
where $r(i)$ is the sum of the elements of the $i$th row of $A$.

Rather than using the Bregman factors $a!^{1/a}$, we will use the
following factors, defined recursively 
in the following fashion:
\begin{equation}
g(1) := e,\ \ 
  g(a + 1) = g(a) + 1 + \frac{1}{2 g(a)} + 
  \frac{0.6}{g(a)^2} \ \ \forall a \geq 1.
\end{equation}
\begin{lem}
For all 0-1 matrices $A$, with row sums $\bf r$ 
\begin{equation}
\label{EQN:weakmincs}
\tilde M(A) = \prod_{i=1}^n \frac{g(r(i))}{e} \geq per(A).
\end{equation}
\end{lem}

The factors in (\ref{EQN:weakmincs}) are fairly close to the factors
in Bregman's theorem.  Using Stirling's formula, the ratio of the
two can be shown to converge to 1.

\begin{tabular}{|c|cccccccccc|}
\hline
$a$ & 1 & 2 & 3 & 4 & 5 & 6 & 7 & 8 & 9 & 10 \\
\hline
$a^{1/a}$ & 1 & 1.41 & 1.82 & 2.21 & 2.61 & 2.99 & 3.38 & 3.76 & 4.15 & 4.53\\
$g(a) / e$ & 1 & 1.47 & 1.89 & 2.31 & 2.71 & 3.11 & 3.50 & 3.89 & 4.27 & 4.66\\
%$g(a) / e$ 
%  & 1 & 1.47 & 1.91 & 2.33 & 2.74 & 3.14 & 3.54 & 3.93 & 4.32 & 4.71 \\
\hline
\end{tabular}

\paragraph{Van der Waerden's Inequality.}
To lower bound the number of perfect matchings,
we will use Van der Waerden's Inequality~\cite{vanlintw92}, 
which says that any doubly stochastic $n$ by $n$ matrix 
has permanent at least $n! / n^n$.

Suppose that our graph is regular, or equivalently, that all the
row and column sums of the matrix $A$ are $\Delta$.  Then dividing
each row by $\Delta$ leaves us with a new doubly stochastic 
matrix $A'$ with
$per(A') = \Delta^{-n} per(A)$.  Applying Van der Waerden's
inequality gives us $per(A) \geq \Delta^n n! / n^n.$
Using Stirling's
Inequality we have that $n! > n^n e^{-n} \sqrt{2 \pi n}$, 
so our lower bound
becomes $(\Delta / e)^n \sqrt {2 \pi n}$.  This differs from
the upper bound by a factor that is polynomial when $\Delta$ is
of the same order as $n$.

\section{The algorithm.}

The basic form of the algorithm is not new, it is simply
acceptance rejection.  What is new is the form of $\tilde M(A)$ that
allows us to use the selfreducibility of the permanent problem together
with the acceptance/rejection protocol.  We will choose
a number uniformly between 1 and $\tilde M(A)$ inclusive.  If this number
corresponds to a valid perfect matching, we accept,
otherwise we reject and must begin again.  We will build up the
perfect matching (permutation) corresponding to our random number 
one entry at a time.  That is, we move through the columns
one at a time, choosing a different row to go with each column.

Let $\sigma$ denote our permutation.
When choosing row $\sigma(i)$ to go with column $i$, update the
matrix by zeroing out column $i$ and row $\sigma(i)$, except for the entry
$A(\sigma(i),i)$ that remains 1.  This leaves a smaller version
of our original problem.  Let $f(A,i,j)$ denote this new matrix.

\begin{figure}
\label{FIG:example}
\begin{center}
\begin{tabular}{ccccc}
$A$ & $f(A,1,1)$ & $f(A,2,1)$ & $f(A,3,1)$ \\
$\left[ \begin{array}{cccc}
 1 & 0 & 1 & 0 \\
 1 & 1 & 0 & 1 \\
 1 & 1 & 1 & 1 \\
 0 & 0 & 1 & 0 \end{array} \right]$ 
&
$\left[\begin{array}{cccc}
 {\bf 1} & 0 & 0 & 0 \\
 0 & 1 & 0 & 1 \\
 0 & 1 & 1 & 1 \\
 0 & 0 & 1 & 0 \end{array} \right]$ 
&
$\left[ \begin{array}{cccc}
 0 & 0 & 1 & 0 \\
 {\bf 1} & 0 & 0 & 0 \\
 0 & 1 & 1 & 1 \\
 0 & 0 & 1 & 0 \end{array} \right]$ 
&
$\left[ \begin{array}{cccc}
 0 & 0 & 1 & 0 \\
 0 & 1 & 0 & 1 \\
 {\bf 1} & 0 & 0 & 0 \\
 0 & 0 & 1 & 0 \end{array} \right]$ \\
\end{tabular}
\end{center}
\caption{Example of reduced matrices}
\end{figure}

Consider the example of Figure~\ref{FIG:example}  For the first column,
we can choose rows 1, 2, or 3 to add to our permutation.  Each choice
leaves a reduced matrix $f(A,1,1), f(A,2,1)$ or $f(A,3,1)$ from which to
choose the remaining permutation.  If our random number from $1$ to 
$M(A)$ lies in $\{1,\ldots,\tilde M(f(A,1,1))\}$, we set $\sigma(1) = 1$.  If
the number lies in 
$\{\tilde M(f(A,1,1)) + 1,\ldots, \tilde M(f(A,1,1)) + \tilde M(f(A,2,1))\}$ we
set $\sigma(1) = 2$, and so on.  The probability of selecting a particular
permutation will be $1 / \tilde M(A)$, and so conditioned on acceptance the
algorithm chooses a permutation with probability $1 / per(A).$

In practice the user will not choose at the beginning of the algorithm
a number from 1 to $\tilde M(A).$  Instead, the user assigns row $\sigma(j)$ to
column $j$ with probability $\tilde M(f(A,\sigma(j),j)) / \tilde M(A)$ 
at each step.
The probability of ending at a particular permutation $\sigma$ will then be:
\begin{equation}
\label{EQN:factors}
\frac{\tilde M(f(A,\sigma(1),1))}{\tilde M(A)} \cdot 
 \frac{\tilde M(f(f(A,\sigma(2),2)))}{\tilde M(f(A,\sigma(1),1))} \cdots
 \frac{\tilde M(f(f(\cdots f(A,\sigma(n),n)))}
 {\tilde M(f(f(\cdots f(A,\sigma(n - 1),n - 1))} = \frac{1}{\tilde M(A)}.
\end{equation}
Pseudocode for this procedure is given in Figure~\ref{FIG:algorithm}

\begin{figure}[h]
\label{FIG:algorithm}
\begin{center}
\vspace{.1in}
\begin{tabular}{|l|}
\hline
{\bf Generate Random Perfect Matching}\\
\\
{\em Input:}  Square 0-1 $n \times n$ 
 matrix $A$ \\
{\em Output:} Perfect matching $\pi$.
\\
{\bf Repeat} \\
\tab {\bf Let} ACCEPT $\leftarrow$ TRUE, 
 {\bf Let} $\tilde A \leftarrow A$, {\bf Let} $\bf {\tilde r}$ be
the row sums of $\tilde A$ \\
\tab {\bf For} $j$ from 1 to $n$\\
\tab \tab {\bf If} $\tilde r(i) = 0$ for any $i$ or
there exists $i \neq i'$ such that 
  $\tilde A(i,j) = \tilde A(i',j) = r(i) = r(i') = 1$\\
\tab \tab \tab {\bf Let} ACCEPT $\leftarrow$ FALSE \\
\tab \tab {\bf Else} \\
\tab \tab \tab {\bf Choose} $R$ at random from 
   $\{1,\ldots,n+1\}$ using \\
\tab \tab \tab \tab 
   $P(R = i) = [\tilde M(f(\tilde A,i,j)) / \tilde M(\tilde A)]
   1(\tilde A(i,j) = 1)$ 
   for $1 \leq i \leq n$ \\
\tab \tab \tab \tab $P(R = n + 1) = 1 - \sum_i P(R = i)$ \\
\tab \tab \tab {\bf If} $1 \leq R \leq n$ \\
\tab \tab \tab \tab {\bf Let} $\pi(j) \leftarrow R$ \\
\tab \tab \tab \tab {\bf Let } $\tilde A \leftarrow f(\tilde A,R,j)$,
     {\bf Let} $\bf {\tilde r}$ be the row sums of $\tilde A$ \\ 
\tab \tab \tab {\bf Else} \\
\tab \tab \tab \tab {\bf Let} ACCEPT $\leftarrow$ FALSE \\
{\bf Until} ACCEPT = TRUE \\
\\
\hline
\end{tabular}
\end{center}
\caption{Generate a perfect matching}
\end{figure}

Important note: to have a valid algorithm we need
$\sum_{i:A(i,j) = 1} \tilde M(f(A,i,j)) / \tilde M(A) \leq 1$.  Otherwise the
factors in the product in (\ref{EQN:factors}) will not be probabilities.
In fact, this is the very reason we use $\tilde M(A)$ rather than $M(A).$  
Suppose that $A$ is a 5 by 5 matrix whose row sums are all 4 and first column
is all ones.  Then $\sum_{i:A(i,j) = 1} M(f(A,i,1)) = 54.5,$ but 
$M(A) = 53.1$.  Since $54.5 / 53.1 > 1$, the algorithm would fail at this
point.  However, for this example
$\sum_{i:A(i,j) = 1} \tilde M(f(A,i,1)) = 64.3,$ and $\tilde M(A) = 65.1,$ 
and so our algorithm proceeds without difficulty.  In 
Section~\ref{SEC:analysis} we 
show in general
that $\sum_{i:A(i,j) = 1} \tilde M(f(A,i,1)) \leq M(A)$.

\subsection{Estimating the permanent.}
Jerrum, Valiant, and Vazirani showed how to turn a sampling algorithm
into an approximate counting algorithm for selfreducible
problems~\cite{jerrumvv86}, but we will not need to
employ their method here.

Instead, we simply note that $per(A) / \tilde M(A)$ is exactly the
probability that the algorithm creates a valid permutation on 
any given run through the columns.  $\tilde M(A)$ is of course
easy to compute in $O(n \log |\tilde M(A)|)$ time, we then just keep track
of the number of acceptances over the number of attempts, and
multiplying by $\tilde M(A)$ gives an approximation to $per(A)$.
Standard Chernoff bounds~\cite{chernoff52} 
show that after $\sigma^{-2} \log(1/\delta)$ samples, the
estimate will be correct to within a factor of $1 + \sigma$ with probability
at least $1 - \delta$.

\section{Analyzing the algorithm.}
\label{SEC:analysis}
Because we are just sampling uniformly from the numbers 
$1,\ldots,\tilde M(A),$
if our algorithm succeeds, it reaches any particular
permutation with probability $1 / \tilde M(A)$, and conditional on reaching
a valid permutation (by which we mean 
one that corresponds to a perfect matching), 
the probability of hitting
a particular permutation is just 
$[1 / \tilde M(A)] / [per(A) / \tilde M(A)] = 1 / per(A)$ exactly as desired.

The only remaining question is whether or not each step in
the algorithm can be successfully performed.  
\begin{lem}
\label{LEM:goal}
For all 0-1 $n$ by $n$ matrices $A$, and $j \in \{1,\ldots,n\}$,
\begin{equation}
\label{EQN:goal}
\sum_{i=1}^n A(i,j) \tilde M(f(A,i,j)) \leq \tilde M(A).
\end{equation} 
\end{lem}

Begin the proof of Lemma~\ref{LEM:goal} by noting that $\tilde M(f(A,i,j))$
and $\tilde M(A)$ each contain many of the same (positive) terms that can
be canceled on both sides of the equation.  Without loss of generality,
assume that elements $A(1,j)$ through $A(c_j,j)$ are 1 and the remaining
elements of column $j$ are zero.  Then when we consider $\tilde M(f(A,i,j))$
for $i \leq c_j$, the $i$ term goes from $g(r_i)$ down to $g(1)$,
$c_j - 1$ row sums of $A$ are reduced by 1.  When $i > c_j$, row $i$ is
the same in $f(A,i,j)$ as in $A$, and so the row sum is also the same.
Hence
$$\frac{\tilde M(f(A,i,j))}{\tilde M(A)} = 
  \frac{e}{g(r(k) - 1)} \prod_{k \leq c_j} 
\frac{g(r(k) - 1)}{g(r(k))}.$$
Therefore we have
\begin{equation}
\label{EQN:goal2}
\frac{\sum_i A(i,j) \tilde M(f(A,i,j))}{M(A)} = 
e \left[\prod_{k = 1}^{c_j} 
\frac{g(r(k) - 1)}{g(r(k))} \right]
\left[\sum_{i=1}^{c_j} \frac{1}{g(r(k) - 1)} \right],
\end{equation}
and we are interested in showing that this expression is bounded above
by 1.  We can rewrite the RHS by denoting $g(r(k) - 1) = a_k$, and using
$g(r(k)) = a_k + 1 + 0.5/a_k + 0.2/a_k^2.$  Set
\begin{eqnarray}
\label{EQN:s}
s(a_k) & := & 1 / a_k \\
\label{EQN:p}
p(a_k) & := & 1 / [1 + 1 / a_k + 0.5 / a_k^2 + .6 / a_k^3 \\
h({\bf a}) & := &
  e \left[\prod_{k = 1}^{c_j} p(a_k) \right]
  \left[ \sum_{i=1}^{c_j} s(a_k) \right].
\end{eqnarray}
When ${\bf a}$ is the vector of row sums for the first $c_j$ rows
of $A$, then $h({\bf a})$ equals the RHS of (\ref{EQN:goal2}), and
so our goal is to show that $h({\bf a})$ is bounded above by 1.  

In fact, 
we will show that $h({\bf a}) \leq 1$ for all ${\bf a}$ such that
$g(1) \leq a_k \leq g(n - 1)$ for all $k$.  That is we allow $\bf a$
to vary continuously over the region rather than restricting it to
the finite set of values $\{g(1),\ldots,g(n-1)\}$.  
Since the region is compact
and $h$ continuous, $h$ attains its maximum at a particular 
point ${\bf a}^*$.

Suppose that we evaluate $h$ at a vector ${\bf a}$ so that $a_k < a_\ell$
for some $k \neq \ell.$  Consider the factors and terms of $h$ that
depend on $a_k$ and $a_\ell.$  Given $\bf a$, 
construct a new vector ${\bf \bar a}$
where $\bar a_k = a$ and $\bar a_\ell$ is chosen so that
\begin{equation}
\label{EQN:constraint}
s(a_k) + s(a_\ell) = s(a) + s(\bar a_\ell).
\end{equation}
Let 
$h_{k,\ell}(a) = h(\bar {\bf a}).$  Then we can write
$h_{k,\ell}(a) = C p(a) p(\bar a_\ell),$ where $C$ is a constant.  
Differentiating, we find
$$
h_{k,\ell}'(a) = 
  C[p'(a)p(\bar a_\ell) - p(a) p(\bar a_\ell) s'(a) / s'(\bar a_\ell)].
$$
Since $s'(x) < 0,$ $p(x) > 0$, and $p'(x) > 0$ 
over the region we are interested in, 
$sgn(h'_{k,\ell}(a)) = 
sgn(p'(a)/[p(a) (-s'(a))] - p'(\bar a_\ell) / 
  [p(\bar a_\ell)(-s'(\bar a_\ell))].$  
Note $h_{k,\ell}'(a)$ has a zero at $a = \bar a_\ell$.  Consider the function
$t(x) = p'(x) / [p(x) (-s'(x))].$  Using (\ref{EQN:s}) and (\ref{EQN:p}) and
simplifying, we find that
$t(x) = 1 + [1.3x - .6]/[x^3 + x^2 + .5x + .6].$  Since $x \geq g(1) = e$ the
numerator is positive, and the denominator is growing faster than the 
numerator so $t(x)$ is a decreasing function over our region of interest.
This
means that $h'_{k,\ell}(a)$ 
is 0 exactly when $a = \bar a_\ell$, is negative when 
$a < \bar a_\ell,$ and positive when $a > \bar a_\ell$.

Hence $h_{k,\ell}(a)$ has its unique maximum at $a = \bar a_\ell$.
Since this is true for any $k$ and $\ell$ with $a_k < a_\ell$, 
all the components of
$\bf a^*$ are identical.  

So instead of maximizing $h({\bf a})$, we can maximize 
\begin{equation}
\bar h(a) := e p(a)^{c_j} \frac{c_j}{a}.
\end{equation}

To accomplish this we will need two facts about exponentials that
are easily verified by taking the appropriate derivatives.
\begin{equation}
\label{EQN:inequality1}
(\forall x \geq 0)(1 - x \leq \exp\{-(x + 0.5 x^2)\}).
\end{equation}
\begin{equation}
\label{EQN:inequality2}
(\forall y \in {\bf R})(y \exp\{-y\} \leq \exp\{-1\}).
\end{equation}

Let $\delta(a) := [1/a + .5/a^2 + .6/a^3]/[1 + 1/a + .5/a^2 + .6/a^3]$, 
so
$p(a) = 1 - \delta(a)$.  Using (\ref{EQN:inequality1}) and 
(\ref{EQN:inequality2}), we have that
\begin{eqnarray*}
\bar h(a) & \leq & \frac{e}{a}c_j\exp\{-c_j(\delta(a) + .5\delta(a)^2)\} \\
& \leq & \frac{1}{a (\delta(a) + .5 \delta(a)^2)}
\end{eqnarray*}
We can write $1 / [a(\delta(a) + .5 \delta(a)^2)]$ as $p_1(a) / p_2(a)$,
where $p_1$ and $p_2$ are degree six polynomials positive for positive $a$.  
Furthermore,
$p_2(a) - p_1(a) = .2a^4 - .05a^3 + .35 a^2 + 1.2a - 3.6$ which is positive
for $a \geq g(1)$, and so 
$p_1(a) / p_2(a) = 1 - [p_2(a) - p_1(a)]/p_2(a) \leq 1.$
Hence $\bar h(a) \leq 1$, and Lemma~\ref{LEM:goal} is proved.

\subsection{Running Time.}

In this section we derive an upper bound on the expected running
time of the algorithm that is polynomial under certain conditions.  
Even if the graph does not meet these conditions,
the algorithm can be used to 
generate perfect matchings.   We only lack a priori
bounds on what the running time is.

The time spent inside the repeat loop is
easy to compute.  Computing the ratios 
$\tilde M(f(\tilde A,i,j)) / \tilde M(\tilde A)$ takes time $O(n)$.
Choosing $R$ takes only $O(n)$ time and $O(\log n)$ 
expected random bits.  Marking the permutation and changing row $R$ 
and column $j$ to
0 also takes $O(n)$ time.
The loop is run $n$ times, so altogether
$O(n^2)$ work is needed. 

The only question that remains is how often does the algorithm accept,
that is, how small is $per(A) / \tilde M(A)$?  The expected amount of samples
taken before acceptance 
will be $\tilde M(A) / per(A)$.
We noted in section~\ref{SEC:I} that
Van der Waerden's inequality~\cite{vanlintw92} together
with Stirling's formula may be used to
show that $per(A) \geq (\sqrt{2 \pi n}) (\Delta / e)^n$.

Since $g(a + 1) \geq g(a) + 1$, we know that $g(a) \geq a$, and so
$g(a + 1) - g(a) < 1 + .5/a + .6/a^2.$  Given this fact, it is straightforward
to show inductively that for $a \geq 2$,
$$g(a) \leq a + .5 \ln a + 1.65.$$
With each row sum identically $\Delta$, 
$\tilde M(A) \leq ([\Delta + .5\ln \Delta + 1.65] / e)^n.$  Using 
$1 + x\leq e^x,$
\begin{eqnarray*}
\frac{\tilde M(A)}{per(A)} \leq 
  \frac{1}{\sqrt{2\pi n}}\left(1 + .5 \frac{\ln \Delta}{\Delta} + 
  \frac{1.65}{\Delta} \right)^n \leq 
  \frac{1}{\sqrt{2\pi n}} [\sqrt{\Delta} + 5.3]^{n / \Delta}
\end{eqnarray*}
We have proved the following
\begin{thm}
The expected running time needed to obtain a sample is
\begin{equation}
O\left(n^{1.5} \Delta^{.5 n / \Delta} 5.3^{n / \Delta}\right) 
\end{equation}
In particular, if $\Delta = \gamma n$ for some constant $\gamma$,
then the expected running time is
$$O(n^{1.5 + .5/\gamma}).$$
\end{thm}

\paragraph{Running time using Markov chains.}  Jerrum and 
Sinclair~\cite{jerrums89} showed
that the Markov chain of Broder~\cite{broder86}
requires $O(n^2 |E| F \ln n)$ time to generate an approximate sample,
where $|E|$ is the number of edges in the graph and $F$ is the ratio
of almost matchings (with $n-1$ edges) to perfect matchings.  

One way to count the number of almost matchings is just to choose
which node on the left is unmatched (there are $n$ ways to do this)
and then use our algorithm to complete the almost matching as before.
This shows that the number of almost matchings is at most $n \tilde M(A)$,
hence the ratio of almost matchings to perfect matchings is at most
$n M(r) / per(A)$, and the total time needed by the Broder chain will
be worse than our algorithm by a factor of $O(n^3 \ln n)$.

\subsection{Estimating the Permanent.}

One form of Chernoff's Bound~\cite{chernoff52} is the following:
\begin{thm} {\bf Chernoff's Bound}
Let $X_1,\ldots,X_t$ be ${0,1}$ i.i.d. Bernoulli random variables 
with parameter $p$.  Then for $\sigma < 1$ 
$$P\left(|\sum_i X_i - tp| > \sigma tp\right) < 2e^{-\sigma^2 tp/3}.$$
\end{thm}

In our algorithm, $p = per(A) / \tilde M(r)$.
Hence after $O([\tilde M(r) / per(A)] \log(1/\delta) / \sigma^2)$ steps,
the algorithm will come within $1 + \sigma$ of the true answer 
with probability at least $1 - \delta$.

\begin{thm}
The expected running time needed to obtain a $1 + \sigma$ 
approximation with probability at least $1 - \delta$ is
$$O\left(n^{1.5 + .5/\gamma}
  \log(1/\delta) / \sigma^2 \right).$$
\end{thm}

\section{Nearly Regular Graphs.}
Until now, we have assumed that the graph was regular so that
each node has the same degree.  We now show how to relax this
assumption to deal with graphs that are merely close to being regular.

In obtaining our upper bound, we did not utilize the regularity
assumption at all, so the only time where it comes into play is
in the lower bound.  Van der Waerden's Inequality is easily 
extended to the case where all the row and column sums are 1 or larger.

Suppose that when we normalize the row sums by dividing by
the degree, the column sums are not exactly 1.  For each column
with sum less than 1, we normalize it by dividing by its column sum.
Then all column and row sums are at least 1, and Van der Waerden's
Inequality assures us that the permanent will be at least 
$n! / n^n$.  Taking account of the normalization, the 
original permanent will be at least
$(n!/ n^n) \prod_{c_i < 1} c_i$.

So if one or two node have degrees that are different from the rest
by a factor of 2, the running time of this approach is unchanged.
If all the degrees are in the set $\{\Delta - c,\ldots, \Delta + c\}$
for some constant $c$, then each 
$c_i \geq (\Delta - c)/ (\Delta + c)$, and their product (for
sufficiently large $n$ will be at least $e^{-2 c n / \Delta}$.  When
$n / \Delta$ is a constant, this is also a constant so the running
time of the procedure will be of the same order as before.
A similar argument shows that even if the degrees of all of the 
nodes vary by $O(\ln n)$, the running time of the procedure will
still be polynomial.  

It can be shown that any bipartite perfect matching problem may be 
efficiently reduced
to a perfect matching problem where the degree of the edges is
either 2 or 3.  Using a technique of Broder~\cite{broder86} we
can then make the problem dense, and the degree of each edge will
still only differ by 1, making it nearly regular.  
Hence the nearly regular problems we are
approximating in this section are still $\sharp P$ complete.

It remains to be seen if this approach of sampling from Minc as
an upper bound can be proven efficient in situations where 
the degrees are not close to being regular.  The denseness assumption
might make it possible to obtain a better lower bound than that given
by Van der Waerden.  

Our algorithm runs quickly on any matrix where 
$\tilde M(r) / per(A)$ is polynomial.  The
Jerrum Sinclair results show that the Markov chain approach runs
quickly when the 
ratio of (n-1)-matchings (matchings with only n-1 edges)
to perfect matchings is polynomial, 
and it would of course be useful to be to make 
precise the connection between these two ratios.

\bibliographystyle{plain}
\bibliography{../../references/myrefsabbrev}

\begin{thebibliography}{10}

\bibitem{barvinok1999}
A.~Barvinok.
\newblock Polynomial time algorithms to approximate permanents and mixed
  discriminants within a simply exponential factor.
\newblock {\em Random Structures Algorithms}, 14:29--61, 1999.

\bibitem{bregman73}
L.~M. Bregman.
\newblock Some properties of nonnegative matrices and their permanents.
\newblock {\em Soviet. Math. Dokl.}, 14(4):945--949, 1973.

\bibitem{broder86}
A.Z. Broder.
\newblock How hard is it to marry at random? ({O}n the approximation of the
  permanent).
\newblock In {\em Proc. 18th ACM Sympos. on Theory of Computing}, pages 50--58,
  1986.

\bibitem{chernoff52}
H.~Chernoff.
\newblock A measure of asymptotic efficiency for test of a hypothesis based on
  the sum of observations.
\newblock {\em Ann. of Math. Stat.}, 23:493--509, 1952.

\bibitem{diaconiss91}
P.~Diaconis and D.~Stroock.
\newblock Geometric bounds for eigenvalues of {M}arkov chains.
\newblock {\em Ann. Appl. Probab.}, 1:36--61, 1991.

\bibitem{efronp1999}
B.~Efron and V.~Petrosian.
\newblock Nonparametric methods for doubly truncated data.
\newblock {\em J. Amer. Statist. Assoc.}, 94(447):824--834, 1999.

\bibitem{jerrums89}
M.~Jerrum and A.~Sinclair.
\newblock Approximating the permanent.
\newblock {\em J. Comput.}, 18:1149--1178, 1989.

\bibitem{jerrumsv2001}
M.~Jerrum, A.~Sinclair, and E.~Vigoda.
\newblock A polynomial-time approximation algorithm for the permanent of a
  matrix with non-negative entries.
\newblock In {\em Proc. 33rd ACM Sympos. on Theory of Computing}, pages
  712--721, 2001.

\bibitem{jerrumvv86}
M.~Jerrum, L.~Valiant, and V.~Vazirani.
\newblock Random generation of combinatorial structures from a uniform
  distribution.
\newblock {\em Theoret. Comput. Sci.}, 43:169--188, 1986.

\bibitem{jerrumv96}
M.~Jerrum and U.~Vazirani.
\newblock A mildly exponential approximation algorithm for the permanent.
\newblock {\em Algorithmica}, (4/5):392--401, 1996.

\bibitem{vanlintw92}
J.~H.~Van Lint and R.~M. Wilson.
\newblock {\em A Course in Combinatorics}.
\newblock Cambridge University Press, 1992.

\bibitem{minc1963}
H.~Minc.
\newblock Upper bounds for permanents of (0,1)-matrices.
\newblock {\em Bull. Amer. Math. Soc.}, 69:789--791, 1963.

\bibitem{sinclair92}
A.~Sinclair.
\newblock Improved bounds for mixing rates of {M}arkov chains and
  multicommodity flow.
\newblock {\em Combin. Probab. Comput.}, 1:351--370, 1992.

\bibitem{valiant79}
L.~G. Valiant.
\newblock The complexity of computing the permanent.
\newblock {\em Theoret. Comput. Sci.}, 8:189--201, 1979.

\end{thebibliography}

\end{document}